\documentclass[10pt]{amsart}
\usepackage{amsmath,amscd,amssymb}
\newtheorem{theorem}{Theorem}[section]
\newtheorem{corollary}[theorem]{Corollary}
\newtheorem{proposition}[theorem]{Proposition}

\newtheorem{lemma}[theorem]{Lemma}
\theoremstyle{remark}
\newtheorem{remark}[theorem]{Remark}
\newtheorem{definition}[theorem]{Definition}
\newtheorem{remarks}[theorem]{Remarks}
\newtheorem{example}[theorem]{Example}
\newtheorem{examples}[theorem]{Examples}
\newtheorem{exercise}[theorem]{Exercise}
\newcommand\A{\mathcal{A}}
\newcommand\be{\begin{equation}\label}
\newcommand\ee{\end{equation}}

\newcommand{\V}{\mathbb{V}}
\newcommand{\TM}{\mathbb{T}M}

\renewcommand{\O}{\mathcal{O}}
\newcommand{\Co}{\mathcal{C}}
\newcommand{\T}{\mathbb{T}}

\newcommand{\X}{\mathcal{X}}

\newcommand{\R}{\mathbb{R}}
\newcommand{\C}{\mathbb{C}}

\newcommand{\Z}{\mathbb{Z}}

\newcommand{\pr}{\on{pr}}
\newcommand{\Cl}{{\on{Cl}}}
\newcommand\lie[1]{\mathfrak{#1}}
\renewcommand{\k}{\lie{k}}

\newcommand{\g}{\lie{g}}

\newcommand{\on}{\operatorname}

\newcommand{\Ad}{ \on{Ad} } 
\newcommand{\ad}{ \on{ad} }

\newcommand{\Hom}{ \on{Hom}} 
\renewcommand{\ker}{ \on{ker}} 
 \newcommand{\Spin}{ \on{Spin}}
\newcommand{\SU}{ \on{SU}} 
\newcommand{\SO}{ \on{SO}}
\newcommand{\GL}{ \on{GL}}

\newcommand{\Mult}{ \on{Mult}}

\newcommand\dirac{/\kern-1.2ex\partial} % Dirac operator
\newcommand\qu{/\kern-.7ex/} % Categorical quotients
\newcommand{\lra}{\longrightarrow}
\newcommand{\hra}{\hookrightarrow}

\renewcommand{\d}{{\mbox{d}}}
\newcommand{\ol}{\overline}

\newcommand\Phinv{\Phi^{-1}}

\newcommand\Sig{\Sigma}
\newcommand\sig{\sigma}

\newcommand\Om{\Omega}
\newcommand\om{\omega}

\newcommand{\f}{\frac}

\renewcommand{\l}{\langle}
\renewcommand{\r}{\rangle}
\newcommand{\hh}{{\textstyle \f{1}{2}}}
\newcommand{\ti}{\tilde}

\newcommand\pt{\on{pt}}

\newcommand\beqn{\begin{equation}}      
\newcommand\eeqn{\end{equation}}      
\newcommand{\ca}{\mathcal}
\newcommand{\wh}{\widehat}
\newcommand{\wt}{\widetilde}
\newcommand{\mf}{\mathfrak}
\newcommand{\beq}{\begin{eqnarray*}}
\newcommand{\eeq}{\end{eqnarray*}}

\newcommand{\Cour}[1]      {[\![#1]\!]}
\newcommand{\ts}{\textstyle}
\newcommand{\Lag}{\on{Lag}}
\newcommand{\Fact}{{\tt FACT}}
\begin{document}

\title[]{Lectures on Pure Spinors and moment maps}

\author{E. Meinrenken}
\address{University of Toronto, Department of Mathematics,
40 St George Street, Toronto, Ontario M4S2E4, Canada }
\email{mein@math.toronto.edu}

\date{\today}

\maketitle 

\tableofcontents

\vskip.5in
\section{Introduction}
This article is an expanded version of notes for my lectures at the
summer school on `Poisson geometry in mathematics and physics' at Keio
University, Yokohama, June 5--9 2006. The plan of these lectures was
to give an elementary introduction to the theory of Dirac structures,
with applications to Lie group valued moment maps. Special emphasis
was given to the \emph{pure spinor} approach to Dirac structures,
developed in Alekseev-Xu \cite{al:der} and Gualtieri \cite{gua:ge}.
(See \cite{bur:int,bur:poi,cou:di} for the more standard approach.)
The connection to moment maps was made in the work of Bursztyn-Crainic
\cite{bur:di}. Parts of these lecture notes are based on a
forthcoming joint paper \cite{al:pur} with Anton Alekseev and Henrique
Bursztyn.

I would like to thank the organizers of the school, Yoshi Maeda and
Guiseppe Dito, for the opportunity to deliver these lectures, and for
a greatly enjoyable meeting. I also thank Yvette Kosmann-Schwarzbach
and the referee for a number of helpful comments.

\vskip.5in
\section{Volume forms on conjugacy classes}\label{sec:intro}
We will begin with the following \Fact, which at first sight may seem quite
unrelated to the theme of these lectures:
\begin{quote}
  \Fact. Let $G$ be a simply connected semi-simple real Lie group. Then every
  conjugacy class in $G$ carries a canonical invariant volume form.
\end{quote}
By definition, a conjugacy class $\Co$ is an orbit for the
conjugation action,
\[ \Ad\colon G\to \on{Diff}(G),\ \Ad(g).a=gag^{-1}.\]
It is thus a smooth $\Ad$-invariant submanifold of $G$. By the
existence of a `canonical' volume form, we mean that there exists an
explicit construction, not depending on any further choices.

More generally, the above {\Fact} holds for simply connected Lie
groups with a bi-invariant pseudo-Riemannian metric. In the
semi-simple case, such a metric is provided by the Killing form.  The
assumption that $G$ is simply connected may be relaxed as well -- the precise
condition will be given below.  Without any assumption on $\pi_1(G)$,
the conjugacy classes may be non-orientable, but they still carry
canonical invariant measures.
\begin{exercise}
\begin{enumerate}
\item Show that $\SO(3)$ has a conjugacy class diffeomorphic to $\R
      P(2)$. This is the simplest example of a non-orientable
      conjugacy class. 
\item Let $G$ be the conformal group of the real line $\R$ (the group
      generated by dilations and translations). Show that
      $G$ does not admit a bi-invariant pseudo-Riemannian metric, and
      that $G$ has conjugacy classes not admitting invariant measures. 
\end{enumerate}
\end{exercise}
The above {\Fact}  does not appear to be well-known.  Indeed, it is not
obvious how to use the pseudo-Riemannian metric to produce a measure
on $\Co$, since the restriction of this metric to $\Co$ may be
degenerate or even zero.
Recall on the other hand that any co-adjoint orbit $\ca{O}\subset \g^*$
carries a canonical volume form, the \emph{Liouville form} for the
Kirillov-Kostant-Souriau symplectic form $\omega$ on $\O$:
\begin{equation}\label{eq:KKS}
 \omega(\xi_1^\sharp,\xi_2^\sharp)\big|_\mu=\l\mu,[\xi_1,\xi_2]\r,\ \
 \ \mu\in\O.
\end{equation}
Here $\xi^\sharp\in\mf{X}(\O)$ denotes the vector field generated by
$\xi\in\g$ under the co-adjoint action. Letting $n=\hh \dim\O$, the
Liouville form is $\f{1}{n!}\om^n$, or equivalently the top degree
part of the differential form $\exp\omega$. One is tempted to try something similar for conjugacy classes. Unfortunately, conjugacy
classes need not admit symplectic forms, in general:

\begin{exercise}
  Show that the group $\Spin(5)$ (the connected double cover of
  $\SO(5)$) has a conjugacy class isomorphic to $S^4$. The 4-sphere
  does not admit an almost complex structure, hence also no
  non-degenerate 2-form.
\end{exercise}

Nevertheless, our construction of the volume form on $\Co$ will be
similar to that of the Liouville form on coadjoint orbits. Let 
\[ B\colon\g\times\g\to \R\]
be the $\Ad$-invariant inner product on $\g$, corresponding to
the bi-invariant pseudo-Riemannian metric on $G$. There is an $\Ad$-invariant 2-form
$\omega\in\Omega^2(\Co)$, given by the formula
\begin{equation}\label{eq:GHJW} 
\omega(\xi_1^\sharp,\xi_2^\sharp)\big|_g=B\Big(\textstyle{\f{\Ad_g-\Ad_{g^{-1}}}{2}}\xi_1,\xi_2\Big),\
\ \ g\in\Co.\end{equation}
This 2-form was introduced by Guruprasad-Huebschmann-Jeffrey-Weinstein
in their paper \cite{gu:gr} on moduli spaces of flat connections, and
plays a key role in the theory of group valued moment maps
\cite{al:mom}. Its similarity to the KKS formula \eqref{eq:KKS} becomes evident if we
use $B$ to identify $\g^*$ with $\g$: The KKS 2-form is defined by the 
skew-adjoint operator $\ad_\mu=[\mu,\cdot]$, while the GHJW 2-form is
defined by the skew-adjoint operator $\hh(\Ad_g-\Ad_{g^{-1}})$. 
An important difference is that the GHJW 2-form may well be
degenerate. It can even be zero: 

\begin{exercise}\label{ex:square}
Show that the GHJW 2-form vanishes on the conjugacy class $\Co$ if and only if
the elements of $\Co$ square to elements of the center of $G$. For $G=\SU(2)$,
there is one such conjugacy class (besides the central elements themselves):
$\Co=\{A\in \SU(2)|\ \on{tr}(A)=0\}$.
\end{exercise}

To proceed we need a certain differential form on the group $G$. Let
$\theta^L,\theta^R\in\Omega^1(G)\otimes \g$ denote the left-, 
right-invariant Maurer-Cartan forms. Thus $\theta^L=g^{-1}\d g$ and
$\theta^R=\d g g^{-1}$ in matrix representations of $G$. 

\begin{theorem}\label{th:1}
Suppose $G$ is a simply connected Lie group with a bi-invariant pseudo
Riemannian metric, corresponding to the scalar product $B$ on
$\g$. Then there is a well-defined smooth, $\Ad$-invariant
differential form $\psi\in\Omega(G)$ such that
\begin{equation}\label{eq:psi}
 \psi_g={\det}^{1/2}\big(\textstyle{\f{\Ad_g+1}{2}}\big)\
\exp\Big(\f{1}{4}B\big(\f{1-\Ad_g}{1+\Ad_g}\theta^L,\theta^L\big)\Big)
\end{equation}
at elements $g\in G$ such that $\Ad_g+1$ is invertible. 
\end{theorem}
Note that the 2-form in the exponential becomes singular at points
where $\Ad_g+1$ fails to be invertible. The Theorem ensures that these
singularities are compensated by the zeroes of the
determinant factor. We can now write down our formula for the volume
form on conjugacy classes. 

\begin{theorem}\label{th:2}
With the assumptions of Theorem \ref{th:1}, the top degree part of the
differential form 
\begin{equation}\label{eq:formula}
 e^\omega\iota_\Co^*\psi 
\end{equation}
defines a volume form on $\Co$. Here $\omega\in\Om^2(\Co)$ is the GHJW
2-form on $\Co$, and $\iota_\Co\colon\Co\hra G$ denotes the inclusion.
\end{theorem}
Since $\psi$ is an even form, the Theorem says in particular that
$\dim\Co$ is even.  Although Formula \eqref{eq:formula} is very
explicit, it is not very easy to evaluate in practice. In particular, it
is a non-trivial task to work out the top degree part `by hand', and
to verify that it is indeed non-vanishing! Also, the complicated
formula for $\psi$ may seem rather mysterious at this point.

What I would like to explain, in the first part of these lectures, is
that the differential form $\psi$ is a \emph{pure spinor} on $G$, and 
that Theorem \ref{th:2} may be understood as the non-degeneracy of a
pairing between two pure spinors, $e^{-\omega}$ and
$\iota_\Co^*\psi$ on $\Co$. 

\begin{remark}
  For the case of a \emph{compact} Lie group, with $B$ positive definite,
  Theorem \ref{th:2} was first proved in \cite{al:du}, using a
  cumbersome evaluation of the top degree part of \eqref{eq:formula}.
  The general case was obtained in \cite{al:pur}.
\end{remark}

The volume forms on conjugacy classes are not only similar to the 
Liouville volume form on coadjoint orbits, but are actually 
generalizations of the latter:
\begin{exercise}
  Let $K$ be any Lie group. The semi-direct product $G=\k^*\rtimes
  K$ (where $K$ acts on $\k^*$ by the co-adjoint action) carries a
  bi-invariant pseudo-Riemannian metric, with associated bilinear form $B$ on
  $\g=\k^*\rtimes\k$ given by the pairing between $\k$ and $\k^*$.
  Show that the inclusion $\k^*\hra G$ restricts to a diffeomorphism
  from any $K$-coadjoint orbit $\O$ onto a $G$-conjugacy class $\Co$. 
  Furthermore, the GHJW 2-form on $\Co$ equals the KKS 2-form on $\O$,
  and the volume form on $\Co$ constructed above is just the ordinary
  Liouville form on $\O$. 
\end{exercise}

\vskip.5in
\section{Clifford algebras and spinors}
This Section summarizes a number of standard facts about Clifford
algebras and spinors. Further details may be found in the classic
monograph \cite{ch:al1}.
\subsection{The Clifford algebra}
Let $W$ be a vector space, equipped with a non-degenerate symmetric
bilinear form $\l\cdot,\cdot\r$. Let $m=\dim W$. 
The \emph{Clifford algebra} over $(W,\l\cdot,\cdot\r)$ is the
associative unital algebra, linearly generated by the elements 
$w\in W$ subject to relations
\[ w_1 w_2 + w_2 w_1 =\l w_1,w_2\r\, 1,\ \ \ w_i\in W.\]
Elements of the Clifford algebra $\Cl(W)$ may be written as linear
combinations of products of elements $w_i\in W$. There is a canonical
filtration,
\[ \Cl(W)=\Cl^{(m)}(W)\supset \cdots \Cl^{(1)}(W)\supset \Cl^{(0)}(W)=\R\]
with $\Cl^{(k)}(W)$ the subspace spanned by products of $\le k$
generators.  The associated graded algebra $\on{gr}(\Cl(W))$ is the
exterior algebra $\wedge(W)$.

The Clifford algebra has a $\Z_2$-grading compatible with the algebra
structure, in such a way that the generators $w\in W$ are odd. With
the usual sign conventions for $\Z_2$-graded (`super') algebras, the
defining relations may be written $[w_1,w_2]=\l w_1,w_2\r\,1$ where
$[\cdot,\cdot]$ denotes the super-commutator. 
A \emph{module over the Clifford algebra}
$\Cl(W)$ is a vector space $\ca{S}$ together with an algebra
homomorphism
\[ \varrho\colon \Cl(W)\to \on{End}(\ca{S}).\] 
Equivalently, the module structure is described by a linear
map $\varrho\colon W\to \on{End}(\ca{S})$ such that 
\[ \varrho(w)\varrho(w')+\varrho(w')\varrho(w)=\l w,w'\r\ 1\]
for all $w,w'\in W$. A Clifford module $\ca{S}$ is called a
\emph{spinor module} if it is irreducible, i.e. if there are no
non-trivial sub-modules. 

\subsection{The Pin group}
Let $\Pi\colon \Cl(W)\to \Cl(W)$ be the parity automorphism of
$\Cl(W)$, equal to $+1$ on the even part and to $-1$ on the odd
part. The \emph{Clifford group} $\Gamma(W)$ is the subgroup of the
group $\Cl(W)^\times$ of invertible elements, consisting of all $x$
such that the transformation 
\begin{equation}\label{eq:transf} y\mapsto\Pi(x)yx^{-1}\end{equation}
of $\Cl(W)$ preserves the subspace $W\subset \Cl(W)$. Let $A_x\in
\GL(W)$ denote the induced transformation of $W$.

\begin{proposition} The homomorphism $A\colon \Gamma(W)\to \GL(W)$ has kernel $\R^\times$ and range $\on{O}(W)$. Thus, the Clifford group fits into an exact sequence, 
\[ 1\lra \R^\times\lra \Gamma(W)\lra \on{O}(W)\lra 1.\]
\end{proposition}
\begin{exercise}
  Show that any $w\in W$ with $B(w,w)\not=0$ lies in $\Gamma(W)$, with
  $A_w$ the reflection defined by $w$. Since any element of
  $\on{O}(W)$ may be written as a product of reflections
  (E.Cartan-Dieudonn\'e theorem), conclude that any element in
  $\Gamma(W)$ is a product $g=w_1\cdots w_k$ with $B(w_i,w_i)\not=0$.
  Use this to prove the above Proposition.
\end{exercise}
Let $x\mapsto x^\top$ denote the canonical anti-homomorphism of $\Cl(W)$,
i.e. $(w_1\cdots w_k)^\top=w_k\cdots w_1$ for $w_i\in W$. Then $g^\top g\in\R^\times$ for all
$g\in\Gamma(W)$. Letting
\[ \on{Pin}(W)=\{g\in \Gamma(W)|\ g^\top g=\pm 1\}\]
one obtains an exact sequence, 
\[  1\lra \Z_2\lra \on{Pin}(W)\lra \on{O}(W)\lra 1.\]

Thus $\on{Pin}(W)$ is a double cover of $\on{O}(W)$. Its restriction
to $\SO(W)$ is denoted $\on{Spin}(W)$.

\subsection{Lagrangian subspaces}
For any subspace $E\subset W$, we denote by $E^\perp$ the space of
vectors orthogonal to $E$. The subspace $E$ is called
\emph{Lagrangian} if $E=E^\perp$.  Let $\on{Lag}(W)$ denote the
\emph{Lagrangian Grassmannian}, i.e. the set of Lagrangian subspaces.
If $\on{Lag}(W)\not=\emptyset$, the bilinear form $\l\cdot,\cdot\r$ is
called \emph{split}. Since we are working over $\R$, the
non-degenerate symmetric bilinear forms are classified by their
signature, and $\l\cdot,\cdot\r$ is split if and only if the signature
is $(n,n)$. That is, $(W,\l\cdot,\cdot\r)$ is isometric to
$\R^{n,n}$, the vector space $\R^{2n}$ with the bilinear form
\[ \l e_i,e_j\r=\pm \delta_{ij},\ \ i,j=1,\ldots,2n\]
with a $+$ sign for $i=j\le n$ and a $-$ sign for $i=j>n$. 
\begin{exercise}\label{ex:lag}
For any invertible matrix $A\in \GL(n)$ let 
\[ E_A=\{(Av,v)|\ v\in \R^n\}\subset \R^{n,n}.\]
\begin{enumerate}
\item Show that $E_A$ is Lagrangian if and only if $A\in \on{O}(n)$. 
\item Show that every Lagrangian subspace $E$ is of the form $E_A$ for 
a unique $A\in\on{O}(n)$. 
%\item Show that the group $\on{O}(n,n)$, in fact already its maximal
%compact subgroup $\on{O}(n)\times\on{O}(n)$, acts transitively on the
%set of Lagrangian subspaces. The stabilizer of the graph of
%$I\in\on{O}(n)$ under this action is the diagonal $\on{O}(n)$.
\end{enumerate}
\end{exercise}
We will assume for the rest of this Section that the bilinear form
$\l\cdot,\cdot\r$ on $W$ is split, of signature $(n,n)$.  The exercise
shows that the Lagrangian Grassmannian $\on{Lag}(W)$ is diffeomorphic
to $\on{O}(n)$. In particular, it is a manifold of dimension
$n(n-1)/2$, with \emph{two connected components}. 
\begin{exercise}
Show that $E,F\in \Lag(W)$ are in the same component of $\Lag(W)$ if
and only if $n+\dim(E\cap F)$ is even. 
\end{exercise}
The orthogonal group
$\on{O}(W)\cong \on{O}(n,n)$ acts transitively on $\Lag(W)$, as does
its maximal compact subgroup $\on{O}(n)\times\on{O}(n)$. (By Exercise
\ref{ex:lag}, already the subgroup $\on{O}(n,0)$ acts transitively.)
\begin{remark}
  Compare with the situation in symplectic geometry: If $(Z,\omega)$
  is a \emph{symplectic} vector space (thus $Z\cong \R^{2n}$ with the
  standard symplectic form $\omega$), a subspace $E$ is called
  Lagrangian if it coincides with its $\omega$-orthogonal space
  $E^\omega$. The symplectic group $\on{Sp}(Z,\omega)$ acts
  transitively on the set $\on{Lag}(Z)$ of Lagrangian subspaces, as
  does its maximal compact subgroup $\on{U}(n)$, and
  $\on{Lag}(Z)=\on{U}(n)/\on{O}(n)$. Thus $\on{Lag}(Z)$ is connected
  and has dimension $n(n+1)/2$. 
\end{remark}

For any pair of transverse Lagrangian subspaces $E,F$, the pairing
$\l\cdot,\cdot\r$ defines an isomorphism $F\cong E^*$. Equivalently,
one obtains an isometric isomorphism 
\[ W\cong E\oplus E^*\]
where the bilinear form on the right hand side is defined by extension
of the pairing between $E$ and $E^*$. 
\begin{exercise}
  Show that for any given $E\in \on{Lag}(W)$, the open subset
  $\{F\in\on{Lag}(W)|\ E\cap F=0\}$ is (canonically) an affine space,
  with $\wedge^2 E$ its space of motions. Show that the closure of
  this subset is a connected component of $\on{Lag}(W)$. Which of the
  two components is it?
\end{exercise}

The sub-algebra of $\Cl(W)$ generated by a Lagrangian subspace $E\in
\Lag(W)$ 
is just the exterior algebra $\wedge E$. Given a Lagrangian
subspace $F$ transverse to $E$, and using the commutator relations to
`write elements of $F$ to the left', we see that
\begin{equation}\label{eq:tensor}
 \Cl(W)=\wedge(F)\otimes\wedge(E),\end{equation}
thus $\Cl(W)=\wedge(W)$ as a $\Z_2$-graded vector space, and also as a
filtered vector space (but not as an algebra). If $w\in W$, the
isomorphism intertwines the operator $[w,\cdot]$ (graded commutator)
on $\Cl(W)$ with the contraction operators $\iota(w)$ on $\wedge(W)$. 

\begin{lemma}
\begin{enumerate}
\item The Clifford algebra $\Cl(W)$ has no non-trivial two-sided
      ideals. 
\item For $E\in\Lag(W)$, the left-ideal $\Cl(W)E$ is maximal. 
\end{enumerate}
\end{lemma}
\begin{proof}
We use the following simple fact: If a non-zero subspace of  
an exterior algebra $\wedge(S)$ is stable under all contraction 
operators $\iota(u),\ u\in S^*$, then the subspace contains the scalars. 

  a) Suppose $\ca{I}$ is a proper 2-sided ideal in $\Cl(W)$. 
Then $\ca{I}$ is invariant under all $[w,\cdot]$ with $w\in W$. 
The above isomorphism \eqref{eq:tensor} takes scalars to scalars, 
and intertwines $[w,\cdot]$ with contractions. Hence $\ca{I}=0$.

b) Let $\ca{I}$ be a proper left-ideal containing $\Cl(W)E$.  
By the
isomorphism \eqref{eq:tensor}, we have a direct sum decomposition
\[\Cl(W)=\wedge(F)\oplus \Cl(W)E.\]
On $\wedge(F)$, the operators $[w,\cdot]$ for $w\in E\cong F^*$
coincide with the contractions $\iota(w)$. Since $\ca{I}\cap
\wedge(F)$ is stable under these operators, it must be zero (or else
it would contain the scalars). Thus $\ca{I}=\Cl(W)E$.
\end{proof}

\begin{corollary}\label{cor:injective}
For any non-zero Clifford module $\ca{S}$ over $\Cl(W)$, 
the action map $\varrho\colon \Cl(W)\to
\on{End}(\ca{S})$ is injective. 
\end{corollary}
\begin{proof}
The kernel of the map $\varrho$ is a 2-sided ideal in $\Cl(W)$, hence
it must be zero.  
\end{proof}

\subsection{The spinor module}
A Clifford module $\ca{S}$ over $\Cl(W)$ is called a
\emph{spinor module} if it is irreducible, i.e. if there are no
non-trivial sub-modules. 
\begin{example}
  If $E\in\on{Lag}(W)$ is Lagrangian, the quotient
  $\ca{S}:=\Cl(W)/\Cl(W)E$ is a spinor module. The irreducibility is
  immediate from the fact that $\Cl(W)E$ is a maximal left-ideal.
\end{example}
We will see below that all spinor modules over $\Cl(W)$ are isomorphic. 
\begin{proposition}
Let $\ca{S}$ be a spinor module, and $E\in\Lag(W)$. Then the 
subspace 
\[ \ca{S}^E=\{\phi\in\ca{S}|\, \varrho(w)\phi=0\  \forall w\in E\}\] 
of elements fixed by $E$ is 1-dimensional. 
\end{proposition}
\begin{proof}
Let $F$ be a complementary Lagrangian
subspace, and choose bases $e_1,\ldots,e_n$ of $E$ and $f^1,\ldots
f^n$ of $F$ with $B(e_i,f^j)=\delta_i^j$. Define $p\in \Cl(W)$ as a product
\[ p=\prod_{i=1}^n e_i f^i=\prod_{i=1}^n (1-f^i e_i)\]
One easily verifies that $p$ has the following properties: 
\begin{equation}\label{eq:pprop}
p^2=p,\ \ Ep=0,\ \ pF=0,\ \ p-1\in \Cl(W)E.
\end{equation}
Since $p^2=p$ the operator $\varrho(p)\in
\on{End}(\ca{S})$ is a projection operator. By $Ep=0$ its
range lies in $\ca{S}^E$, and by $p-1\in \Cl(W)E$ it acts as the
identity on $\ca{S}^E$. Hence
\[ \ca{S}^E=\varrho(p)\,\ca{S}.\]
Since $\varrho$ is injective, we have $\varrho(p)\not=0$, hence
$\ca{S}^E\not=0$. Pick a non-zero element $\phi\in \ca{S}^E$.  Then
$\ca{S}=\varrho(\Cl(W))\phi
=\varrho(\wedge(F))\phi$ by irreducibility, and since 
the left ideal $\Cl(W)E$ acts
trivially on $\phi$. Since $pF=0$, and hence $p\wedge\!(F)=\R p$ 
we obtain
\begin{equation}\label{eq:fix}
 \ca{S}^E=\varrho(p)\,\ca{S}=\varrho(p)\,\varrho(\wedge(F))\,\phi=
\R \varrho(p)\phi=\R\phi.\end{equation}
Equation \eqref{eq:fix} proves that $\ca{S}^E$ is 1-dimensional. 
\end{proof}
The kernel and range of any homomorphism of Clifford modules are
sub-modules. Hence, any non-zero homomorphism of spinor modules is an 
isomorphism. In particular, this applies to the action map
\begin{equation}\label{eq:action}
 \Cl(W)/\Cl(W)E\otimes \ca{S}^E \to \ca{S},\ x\otimes\phi\mapsto
 \varrho(x)\,\phi.
\end{equation}
for any spinor module $\ca{S}$. Thus:
\begin{corollary}
  Any two spinor modules $\ca{S}_1,\ca{S}_2$ over $\Cl(W)$ are
  isomorphic. Furthermore, the isomorphism is unique up to non-zero
  scalar, i.e. the space $\Hom_{\Cl(W)}(\ca{S}_1,\ca{S}_2)$ is
  1-dimensional.
\end{corollary}
As a consequence, the projectivization $\mathbb{P}(\ca{S})$ of the
spinor module is canonically defined (i.e. up to a \emph{unique}
isomorphism). In other words, the Clifford algebra $\Cl(W)$ has a
unique irreducible \emph{projective} module.  The map taking $E$ to
$\ca{S}^E$ defines a canonical equivariant embedding
\[ \on{Lag}(W)\to \mathbb{P}(\ca{S})\] 
as an orbit for the action of $\on{O}(W)$. The image can be
characterized as follows. 
Given $\phi\in\ca{S}$, let
$N_\phi\subset W$ be its `null space', 
\[ N_\phi=\{w\in W|\ \varrho(w)\,\phi=0\}.\]
If $w_1,w_2\in N_\phi$ then 
\[ 0=\varrho(w_1)\varrho(w_2)\phi+\varrho(w_2)\varrho(w_1)\phi=
\varrho([w_1,w_2])\,\phi=B(w_1,w_2)\phi.\]
Hence, if $\phi\not=0$ the subspace $N_\phi$ is isotropic. 
\begin{definition}
  A non-zero spinor $\phi\in\ca{S}$ is called a \emph{pure spinor} if
  $N_\phi$ is Lagrangian. Let $\on{Pure}(\ca{S})$ denote the set of 
  pure spinors  of $\ca{S}$. 
\end{definition}
Note that the pure spinors defining a given Lagrangian subspace $E$
are exactly the non-zero elements of the line $\ca{S}^E$.  We can
summarize the discussion in the following commutative diagram,
equivariant for the action of $\on{Pin}(W)$:
\[ \begin{CD} \on{Pure}(\ca{S}) @>>> \ca{S}^\times\\
@VVV @VVV\\
\Lag(W) @>>> \mathbb{P}(\ca{S})\end{CD}\] 

\begin{exercise}
  Show that for any spinor module, the map $\varrho\colon \Cl(W)\to
  \on{End}(\ca{S})$ is an isomorphism.
\end{exercise}

\begin{exercise}
  Any maximal left ideal $\ca{I}\subset\Cl(W)$ defines a spinor module
  $\Cl(W)/I$. Prove that the set of maximal left ideals is canonically
  isomorphic to $\mathbb{P}(\ca{S})$, and that the inclusion of
  $\Lag(W)$ is just the map $E\to \ca{I}=\Cl(W)E$.
\end{exercise}

\subsection{The Chevalley pairing}\label{subsec:chev}
Let $\ca{S}$ be a spinor module over $\Cl(W)$. Then the dual space $\ca{S}^*$ 
is again a spinor module, with Clifford action given as 
\[ \varrho_{\ca{S^*}}(x)=\varrho_{\ca{S}}(x^\top)^*.\]
We obtain a 1-dimensional line
$K_S=\on{Hom}_{\Cl(W)}(\ca{S}^*,\ca{S})$. The evaluation map defines an
isomorphism of Clifford modules, 
\[ \ca{S}\cong \ca{S}^*\otimes K_S .\]
Tensoring with $\ca{S}$, and composing with the duality pairing
$\ca{S}\otimes\ca{S}^*\to \R$, we obtain a pairing 
\[ \ca{S}\otimes\ca{S}\to K_S,\ \ \phi\otimes\psi\mapsto (\phi,\psi).\]
This pairing is known as the \emph{Chevalley pairing}. By
construction, the Chevalley pairing satisfies 
\[ (\phi,\varrho(x)\psi)=(\varrho(x^\top)\phi,\psi)\]
for all $x\in \Cl(W)$. In particular, 
\begin{equation}\label{eq:equivariance}(g.\phi,g.\psi)=
\pm(\phi,\psi)\end{equation}
for all $g\in \on{Pin}(W)$, where the dot indicates the Clifford action. This just follows since
$g^\top g=\pm 1$ by definition of the Pin group.
\begin{exercise}
  Let $\ca{S}=\Cl(W)/\Cl(W)E$. Show that there is a canonical
  isomorphism $K_S=\wedge^n(E^*)$, where $n=\hh\dim W$. (Hint:
  $\ca{S}^*$ is identified with the submodule of $\Cl(W)$ generated by
  the line $\wedge^n(E)\subset \Cl(W)$.)  Choose a Lagrangian subspace
  complementary to $E$ to identify $\ca{S}=\wedge E^*$. Show that that
  $(\phi,\psi)=(\phi^\top\wedge \psi)_{[n]}$ using the wedge product
  in $\wedge E^*$.
\end{exercise}
\begin{proposition}[Chevalley]\label{th:chevalley}\cite[Theorem III.2.4]{ch:al1}
  If $\phi,\psi$ are pure spinors, then the Lagrangian subspaces
  $N_\phi,N_\psi$ are transverse if and only if $(\phi,\psi)\not=0$.
\end{proposition}
\begin{proof}
  Let $E=N_\phi$. For any Lagrangian complement $F\cong E^*$ to $E$,
  there is a unique isomorphism of spinor modules $\ca{S}\cong \wedge
  E^*$ taking $\phi$ to the pure spinor $1\in\wedge E^*$. In this
  model, $K_S=\wedge^n E^*$, and $(\phi,\psi)=\psi_{[n]}$.  Suppose
  that $(\phi,\psi)\not=0$, thus $\psi_{[n]}\not=0$. 
For $w\in E-\{0\}$ 
we have
\[ (\varrho(w)\psi)_{[n-1]}=\iota_w (\psi_{[n]})\not=0.\]
It follows that $N_\psi\cap
  E=\{0\}$. Conversely, if $N_\psi$ is transverse to $N_\phi$, we may take
  $F=N_\psi$.  Then $\psi\in \wedge^n F^*-\{0\}$, and in particular
  $(\phi,\psi)=\psi_{[n]}\not=0$.
\end{proof}

In particular, pure spinors in $\ca{S}=\Cl(W)/\Cl(W)E$, for any pair of
transverse Lagrangian subspaces define a non-zero element of
$\wedge^nE^*$, i.e. a volume form on $E$.

\vskip.5in
\section{Linear Dirac geometry}
The field of Dirac geometry was initiated by T. Courant in
\cite{cou:di}. One of the original motivations of this theory was to
describe manifolds with `pre-symplectic foliations', arising
for instance as submanifolds of Poisson manifolds.  The term `Dirac
geometry' stems from its relation with the Dirac brackets arising in
this context. One of the key features of Dirac geometry is that it
treats Poisson geometry and pre-symplectic geometry on an equal
footing. More recently, it was observed by Hitchin \cite{hi:gen} that
complex geometry can be understood in this framework as well, leading
to the new field of generalized complex geometry \cite{hi:gen,gua:ge}.

As in Courant's original paper, we will first discuss the linear case.  
\subsection{Linear Dirac structures}
Let $V$ be any vector space, and $\V=V\oplus V^*$ equipped with the
bilinear form $\l\cdot,\cdot\r$ given by the pairing between $V$ and
$V^*$: 
\[ \l v_1\oplus \alpha_1,\ v_2\oplus
\alpha_2\r=\l\alpha_1,v_2\r+\l\alpha_2,v_1\r\]
for $v_i \in V$ and $\alpha_i\in V^*$. Specializing the constructions
from the last section to the case $W=\V$, we note that $\V$ has two
distinguished Lagrangian subspaces, $V$ and $V^*$. We will call the
corresponding spinor modules over $\Cl(\V)$,  
\[ \Cl(\V)/\Cl(\V)V\cong\wedge V^*,\ \ \ \ \Cl(\V)/\Cl(\V)V^*
\cong \wedge V\]
the contravariant and covariant spinor modules, respectively. The star
operator for any volume form on $V$ defines an isomorphism between
these two spinor modules. 
\begin{definition}
A linear Dirac structure on a vector space $V$ is a Lagrangian
subspace $E\in \on{Lag}(\V)$.
\end{definition}
As we have seen, a linear Dirac structure $E$ may be described by a
line $\ca{S}^E$ of pure spinors, using e.g. the covariant or contravariant
spinor module.
\begin{examples}
Consider the contravariant spinor representation $\wedge V^*$. Here are
some examples of pure spinors $\phi$ and associated Lagrangian
subspaces $N_\phi$: 
\begin{enumerate}
\item $\phi=1$ corresponds to $N_\phi=V$.
\item For any 2-form $\omega\in\wedge^2 V^*$, the exponential 
      $\phi=e^{-\omega}$ is a pure spinor, with $N_\phi$ the graph
\[ \on{Gr}_\omega=\{v\oplus \alpha|\ \alpha=\omega(v,\cdot)\}.\]
\item Any volume form $\mu\in \wedge^{\on{top}}V^*\backslash 0$ defines a pure spinor
$\phi=\mu$, with $N_\phi=V^*$. 
\item If $\mu$ is a volume form and $\pi\in\wedge^2 V$, the element 
$\phi=e^{-\iota(\pi)}\mu$ (where $\iota\colon\wedge V\to
\on{End}(\wedge V^*)$ is the algebra homomorphism extending the
contraction operators $v\mapsto \iota(v)$) is a pure spinor, with
$N_\phi$ the graph 
\[ \on{Gr}_\pi=\{v\oplus \alpha|\ v=\pi(\alpha,\cdot)\}.\]
\end{enumerate}
\end{examples}

For any Lagrangian subspace $E\subset \V=V\oplus V^*$, define its \emph{range}
$\on{ran}(E)$ to be the projection onto $V$. One observes that
$S=\on{ran}(E)$ carries a well-defined 2-form,  
\begin{equation}\label{eq:twoform} 
\omega_S(v_1,v_2)=\l\alpha_1,v_2\r=-\l\alpha_2,v_1\r
\end{equation}
where $v_i\oplus \alpha_i\in E$ are lifts of $v_i\in\on{ran}(E)$. The
kernel of this 2-form is $\ker\omega_S=\{v\in V|\ (v,0)\in E\}$.
Conversely, $E$ may be recovered from $S$ together with the 2-form
$\omega_S\in\wedge^2 S^*$, as $E=\{(v,\alpha)|\ v\in S,\ 
\alpha|_S=\omega_S(v,\cdot)\}$.
 
Let
$\on{ann}(S)\subset V^*$ be the annihilator of $S$, and choose a
non-zero element $\mu\in \wedge^{\on{top}}\on{ann}(S)\subset \wedge
V^*$. 
\begin{exercise}\label{ex:standard}
Show that 
\[ \phi=e^{-\omega_S}\mu\in\wedge V^*\]
is a pure spinor with $N_\phi=E$. (Here we have chosen an arbitrary
extension of $\omega_S$ to a 2-form on $V$. Note that the element
$\phi$ does not depend on this choice.) Conversely, show that every
contravariant pure spinor has the form $ \phi=e^{-\omega}\mu$, for uniquely
given $S,\mu\in\wedge^{\on{top}}\on{ann}(S),\omega \in \wedge^2 S^*$.
\end{exercise}   
Put differently, a contravariant pure spinor is equivalent to a Lagrangian
subspace $E$ together with a volume form on $V/\on{ran}(E)$. 

\begin{exercise}
  Work out a similar description for covariant pure spinors.
\end{exercise}

\subsection{Dirac maps}
Let $A\colon V\to V'$ be a linear map. We say that two elements $w=v\oplus \alpha\in
\V$ and $w'=v'\oplus \alpha'\in\V'$ are \emph{$A$-related}, and
write  
\[ w\sim_A w'\] 
if $v'=A(v)$ and $\alpha=A^*(\alpha')$. Then the pull-back map of 
contravariant spinors has the property, 
\[ \varrho(w)(A^*\phi')=A^*(\varrho(w')\phi')\]
for $w\sim_A w'$ and $\phi'\in \wedge (V')^*$. From this, we see that
the pull-back of a contravariant pure spinor $\phi'$ is again a pure
spinor \emph{unless $ A^*\phi'=0$}. Hence, if $F'\subset \V'$ is a
Lagrangian subspace, and $\ca{S}^{F'}\subset \wedge (V')^*$ the pure
spinor line in the contravariant spinor module, then $A^* \ca{S}^{F'}$
is either zero, or is a pure spinor line corresponding to some
Lagrangian subspace $F\subset \V$.
\begin{exercise}
Suppose $\ca{S}^F=A^* \ca{S}^{F'}$. Show that 
\begin{equation}\label{eq:F}
F=\{w\in \V|\ \exists w'\in F'\colon w\sim_A w'\} .\end{equation}
\end{exercise}

Similarly, in the covariant spinor representation we have
\[ \varrho(w')(A_*\chi)=A_*(\varrho(w)\chi)\]
for $w\sim_A w'$ and $\chi\in\wedge V$. Hence, if $E\subset \V$ is a
Lagrangian subspace, and $\ca{S}^E$ is the pure spinor line in the
covariant spinor representation, then $A_*(\ca{S}^E)$ is either zero, or
is the pure spinor line for a Lagrangian subspace $E'\subset \V'$.
\begin{exercise}
Suppose $\ca{S}^{E'}=A_*(\ca{S}^E)$. Show that 
\begin{equation}\label{eq:E}
E'=\{w'\in \V'|\ \exists w\in E\colon w\sim_A w'\}.
\end{equation}
\end{exercise}
\begin{definition}
  Let $V,V'$ be vector spaces with linear Dirac structures $E,E'$. A
  linear map $A\colon V\to V'$ is called a \emph{Dirac map} if the
  spaces $E,E'$ are related by \eqref{eq:E}. It is a 
\emph{strong Dirac map} if
  the induced map $A_*\colon \ca{S}=\wedge(V)\to \ca{S}'=\wedge(V')$ satisfies
  $A(\ca{S}^E)=(\ca{S}')^{E'}$.
\end{definition}
Strong Dirac maps are also called \emph{Dirac realizations} in
the literature.
\begin{exercise}
Show that a Dirac map $A$ is a strong Dirac map if and only if $E\cap (\ker(A)\oplus 0)=0$. 
\end{exercise}

\begin{example}
  Let $\pi\in\wedge^2V,\ \pi'\in\wedge^2 V'$ be 2-forms, with
  $\pi'=A_*\pi$. Then the map $A$ is a strong Dirac map
  relative to $E'=\on{Gr}_{\pi'}$ and $E=\on{Gr}_\pi$.
\end{example}
\begin{example}
Let $E$ be a linear Dirac structure on $V$. Recall that
$S=\on{ran}(E)$ carries a unique 2-form $\omega$, with 
\[E=\{v\oplus \alpha|\ v\in S,\ \alpha|_S=\omega(v,\cdot)\}\] 
View $(S,\omega)$ as a Dirac space, with Dirac structure defined by
the graph $\on{Gr}_\omega\subset S\oplus S^*$. Then the inclusion map 
$\iota_S\colon S\to V$ is a strong Dirac map.
\end{example}

\begin{exercise}
Let $V$ carry the Dirac structure $E$. Then the collapsing map $V\to
\{0\}$ is (trivially) a Dirac map. Show that it is a strong Dirac map
if and only if the 2-form $\omega$ on $S=\on{ran}(E)$ is
non-degenerate, if and only if $E=\on{Gr}_\pi$ for a bi-vector $\pi$.
\end{exercise}

Suppose $E,F\in\Lag(\V)$ are Lagrangian subspaces. Then $E,F$ are
transverse if and only if $\l\psi,\chi\r\not=0$, where $\psi\in\wedge
V^*$ is a contravariant pure spinor defining $F$ and
$\chi\in\wedge(V)$ is a covariant pure spinor defining $E$. (This is
equivalent to Proposition \ref{th:chevalley}.)  The following result
says that Lagrangian complements may be `pulled back' under strong
Dirac maps.

\begin{proposition}\label{lem:trick}
  Suppose $A\colon V\to V'$ is a strong Dirac map relative to
  Lagrangian subspaces $E\subset \V,\ E'\subset \V'$. Let
  $\psi'\in\wedge(V')^*$ be a
  covariant pure spinor, with $N_{\psi'}=F'$ transverse to $E'$. Then 
  $\psi=A^*\psi'$ is non-zero, and $N_\psi=F$ is transverse to $E$. 
  Equivalently, if $\phi$ is a contravariant pure spinor with
  $N_\phi=E$ we have 
\[ (\phi,A^*\psi')\not=0.\]
\end{proposition}
\begin{proof}
Let $\chi\in\wedge V$ be a \emph{covariant} pure spinor defining $E$. Then 
$A_*(\chi)\in\wedge(V')$ is a pure spinor defining $E'$. 
Since  $E',F'$ are transverse,
\[ 0\not=\l\psi',\ A_*\chi\r=\l A^*\psi',\chi\r.\] 
This shows that $\psi=A^*\psi'$ is a pure spinor, with $N_\psi=F$
transverse to $E$. 
\end{proof}
\subsection{The map $\on{O}(V)\to \Lag(\V)$}
Suppose now that $V$ itself carries a non-degenerate symmetric
bilinear form $B$. Let $\ol{V}$ denote the same vector space, but with
the bilinear form $-B$. There is
an isometric isomorphism
\[ \kappa\colon V\oplus \ol{V}\to \V=V\oplus V^*,\ \  v\oplus
w\mapsto (v-w)\oplus \hh B(v+w,\cdot).\]
This identifies $\on{O}(V\oplus \ol{V})\cong
\on{O}(\V)$, and in particular yields an inclusion of the subgroup 
$\on{O}(V)$: 
\[ \on{O}(V)\hra \on{O}(\V),\ A\mapsto A^\kappa=\kappa\circ
\left(\begin{array}{cc}A&0\\0&I\end{array}\right) 
\circ \kappa^{-1}.\]
If we use $B$ to identify $V$ and $V^*$, the matrix on the right is 
easily computed to be, 
\begin{equation}\label{eq:matrix}
 A^\kappa=\left(\begin{array}{cc}(A+I)/2&(A-I)\\(A-I)/4&(A+I)/2\end{array}\right).\end{equation}
Its action on $V$ describes a new Lagrangian subspace,
\[ F=A^\kappa(V).\] 
Let $\Gamma(V)\to \Gamma(\V)$ be the inclusion of Clifford groups
defined by the homomorphism $\Cl(V)\subset \Cl(V\oplus \bar{V})\cong \Cl(\V)$. 
This lifts the map $\on{O}(V)\to \on{O}(\V)$, and restricts to a
homomorphism of $\on{Pin}$ groups. For any lift $\ti{A}\in\Gamma(V)$ of $A$, we obtain a lift $\ti{A}^\kappa
\in \Gamma(\V) $ of $A^\kappa$. The pure spinor $1\in \wedge(V^*)$
defines the Lagrangian subspace $V\subset \V$, hence
\[\psi=\varrho(\ti{A}^\kappa)\,1\] 
represents $F=A^\kappa(V)$. The situation is described in the
following diagram:
\[ \begin{CD} \Gamma(V) @>>> \on{Pure}(\V)\\
@VVV @VVV\\ \on{O}(V) @>>> \on{Lag}(\V)\end{CD}\]
where the lower map is $A\mapsto A^\kappa(v)$ and the upper map is
$\ti{A}\mapsto \varrho(\ti{A}^\kappa)(1)$. Since $\on{Pin}(V)\subset
\Gamma(V)$ is a double cover of $\on{O}(V)$, there is a lift
$\ti{A}\in\on{Pin}(V)$ that is unique up to sign.  One has an
explicit formula for the resulting $\psi$, valid for
$\det(A+I)\not=0$:
\begin{equation}\label{eq:psiform}
 \psi={\det}^{1/2}({\ts{\f{A+I}{2}}})\exp(\textstyle{\f{1}{4}}\sum_i
 (\textstyle{\f{I-A}{I+A}} v_i)\wedge v^i).
\end{equation}
%
%The proof of this formula is based on the following factorization of
%the matrix $A^\kappa$: 
%\[ A^\kappa=\left(\begin{array}{cc} I& 2\f{A-I}{A+I}\\ 0&I\end{array}\right)
%\left(\begin{array}{cc} I& 0\\ \f{A-A^{-1}}{8}&I\end{array}\right)
%\left(\begin{array}{cc} \f{A+I}{2}& 0\\ 0&\f{2}{A^{-1}+I}\end{array}\right)
%\]
See \cite{al:cli} for a proof. 
Here we have used $B$ to identify $V^*\cong V$, and $v_i,v^i$ are
bases with $B(v_i,v^j)=\delta_i^j$. The sign of the square root
depends on the choice of lift $\ti{A}$.

Similarly, the action of $A^\kappa$ on $V^*$ defines a Lagrangian
subspace transverse to $F$, 
\[ E=A^\kappa(V^*).\]
Given an orientation on $V$, the
associated Riemannian volume form $\mu$ is a pure spinor defining
$V^*$, hence 
\[\phi=\varrho(\ti{A}^\kappa)\,\mu\] 
is a pure spinor defining $E$.

\begin{remark}
  If the scalar product $B$ on $V$ is definite, the inclusion
  $\on{O}(V)\to \on{Lag}(V)$ is a bijection: This is exactly the
  diffeomorphism $\on{Lag}(V)\cong \on{O}(n)$ mentioned earlier.
  (This isomorphism is described in the paper   
   \cite{cou:di} under the name `generalized Cayley transform'.) 
  Similarly, the map $\Gamma(V)\to \on{Pure}(\V),\ g\mapsto
  \varrho(g)\,1$ defines a bijection of the set of pure spinors with
  the Clifford group:
\[ \on{Pure}(\V)\cong \Gamma(n):=\Gamma(\R^n).\]
\end{remark}

\vskip.5in
\section{The Cartan-Dirac structure}
\subsection{Almost Dirac structures}
It is straightforward to generalize the above considerations from
vector spaces to vector bundles, and in particular to the tangent
bundle of a manifold. Thus, let 
\[ \T M=TM\oplus T^*M\]
be the generalized tangent bundle, with fiberwise inner product
$\l\cdot,\cdot\r$ given by the pairing of 1-forms with vector fields,
and $\Cl(\T M)$ the corresponding bundle of Clifford algebras.
Covariant spinors are multi-vector fields, $\chi\in
\mf{X}^\bullet(M)=\Gamma(M,\wedge TM)$, while contravariant spinors
are differential forms, $\phi\in \Omega(M)=\Gamma(M,\wedge(T^*M))$.

An \emph{almost Dirac structure} on $M$ is a Lagrangian sub-bundle
$E\subset \T M$. (In Section \ref{sec:dirac}, we will discuss the integrability condition
turning an almost Dirac structure to a Dirac structure.)  A smooth map
between almost Dirac manifolds $f\colon M\to M'$ is called a (strong)
Dirac map if each tangent map $(\d f)_x$ is a (strong) Dirac map.

Any almost Dirac structure may be described (at least locally) by a
contravariant pure spinor $\phi\in\Omega(M)$, or by a covariant pure
spinor $\chi\in\mf{X}(M)$. Our basic examples for vector spaces carry
over to manifolds: Any 2-form on a manifold defines an almost Dirac
structure, as does any bi-vector field. If $E$ is an almost Dirac
structure, described (locally) by a pure spinor $\phi$, and
$\tau\in\Omega^2(M)$ any 2-form, one may define a new almost Dirac
structure $E^\tau$ described locally by pure spinor $e^{-\tau}\phi$.
One calls $E^\tau$ the \emph{gauge transformation} of $E$ by the
2-form $\tau$. For instance, taking
$E=TM$, one obtains the graph of $\tau$:
\[ (TM)^\tau=\on{Gr}_\tau.\] 
\begin{exercise}
In general, show that $E^\tau$ is the image of $E$ under the
automorphism $v\oplus \alpha\mapsto v\oplus (\alpha+\iota_v\tau)$ 
of $\T M$. 
\end{exercise}

Given a pseudo-Riemannian metric on $M$, any section $A$ of the group bundle
$\on{O}(TM)$ defines a pair of transverse Lagrangian sub-bundles
\[ E=A^\kappa(T^*M),\ F=A^\kappa(TM)\] 
of $\T M$. A lift $\ti{A}$ to a section of $\on{Pin}(TM)$ defines pure spinors
$\phi,\psi$ corresponding to $E,F$, where $\phi$ depends on the choice
of an orientation on $M$.

\subsection{The case $M=G$}
Let $G$ be a Lie group. For any $\xi\in\g$, let
$\xi^L,\xi^R\in\mf{X}^1(G)$ the corresponding left-,right-invariant
vector fields. The bundle $\on{GL}(TG)$ has a unique section $A$ with the 
property 
\[ A(\xi^L)=\xi^R.\]
Suppose $G$ carries a bi-invariant pseudo-Riemannian metric,
and let $B$ the corresponding inner product on the Lie algebra $\g$.
Then $A$ is an $\Ad(G)$-invariant section of $\on{O}(TG)\subset \GL(TG)$. Hence, it
determines transverse $\Ad(G)$-invariant Lagrangian sub-bundles $E,F\subset \T G$. 
Recall that
$\theta^L,\theta^R\in\Omega^1(G)\otimes\g$ denote the
Maurer-Cartan forms. 
\begin{proposition}
Define bundle maps $e,f\colon \g\to \T G$ by 
\[ e(\xi)=(\xi^L-\xi^R)\oplus B(\f{\theta^L+\theta^R}{2},\xi)\]
\[ f(\xi)=(\f{\xi^L+\xi^R}{2})\oplus B(\f{\theta^L-\theta^R}{4},\xi)\]
The maps $e,f$ are injective, and have range $E,F$.  
\end{proposition}
\begin{proof}
  Under left-trivialization of the tangent bundle $TG=G\times\g$, the
  section $A$ is just the adjoint action, $g\to \Ad_g$.
  Hence, the section $A^\kappa$ is given by \eqref{eq:matrix}, with
  $A$ replaced by $\Ad_g$. Writing
  elements $\xi_1\oplus\xi_2\in\g\oplus \g\cong T_gG\oplus T_g^*G$ as
  column vectors, we see that 
\[ A^\kappa(\xi_1\oplus
  \xi_2)=f(\xi_1)+e(\xi_2).\]
\end{proof}
The sub-bundle $E$ is called the \emph{Cartan-Dirac structure}. (It
satisfies the integrability condition discussed below.)  Since
$\xi^\sharp=\xi^L-\xi^R$ are the generating vector fields for the
conjugation action, the generalized distribution
$\on{ran}(E)=\on{pr}_{TM}(E)$ is just the distribution tangent to the
conjugacy classes $\Co$ of $G$. Hence, by \eqref{eq:twoform} the conjugacy classes
$\Co\subset G$ acquire $\Ad(G)$-invariant 2-forms $\omega$.
\begin{proposition}\label{prop:ghjw}
The 2-forms $\om$ on conjugacy classes $\Co$ are exactly the 
GHJW 2-forms. 
\end{proposition}
We leave the proof as an exercise. Equivalently, the inclusion maps
\[ \iota_\Co\colon \Co\hra G\]
are strong Dirac maps, relative to the (almost) Dirac
structures given by the GHJW 2-form $\omega$ on $\Co$ and the
Lagrangian sub-bundle $E\subset \T G$. 

Let us now \emph{assume} that the adjoint action $\Ad\colon G\to
\on{O}(\g)$ lifts to a group homomorphism $\wt{\Ad}\colon G\to
\on{Pin}(\g)$ into the Pin group.  (This is automatic if $G$
is simply connected.)  The lift $\wt{\Ad}$ determines lifts
$\ti{A},\ti{A}^\kappa$. Hence it defines an invariant pure spinor
$\psi=\varrho(\ti{A}^\kappa)\,1$ with $N_\psi=F$, and given an
invariant volume form $\mu$ it also defines a pure spinor
$\phi=\varrho(\wt{A}^\kappa)\,\mu$ with $N_\phi=E$.

Consider now a conjugacy class $\Co$, with GHJW 2-form $\omega$. 
By Lemma \ref{lem:trick}, the pull-back
$\iota_\Co^*\psi$ is a pure spinor
defining a Lagrangian sub-bundle transverse to $\on{Gr}_\omega$.
Equivalently, the pairing between the two
pure spinors $e^{-\omega},\iota_\Co^*\psi$ is non-vanishing, that is 
\[
0\not=(e^{-\omega},\iota_\Co^*\psi)=(e^\omega\iota_\Co^*\psi)_{[\on{top}]}\]
is a volume form on $\Co$. We have shown the following more precise
version of the \Fact: 
\begin{theorem}
  Suppose that the adjoint action $\Ad\colon G\to \on{O}(\g)$ lifts to
  a homomorphism $\wt{\Ad}\colon G\to \on{Pin}(\g)$, and let
  $\psi\in\Omega(G)$ be the pure spinor defined by such a lift. Then, for
  any conjugacy class $\Co$ in $G$ the top degree part of
\[  e^\omega\iota_\Co^*\psi\]
defines an invariant volume form on $\Co$. 
\end{theorem} 
The explicit formula
\eqref{eq:psi} for $\psi=\ti{A}^\kappa.1$ is obtained as a special
case from \eqref{eq:psiform}: 
\begin{proposition}
  If $G$ is connected
  and the adjoint action
  $G\to \on{SO}(\g)$ lifts to a group homomorphism $G\to \on{Spin}(\g)$,
  the formula \eqref{eq:psi} defines a pure spinor $\psi$ with
  $N_\psi=F$.
\end{proposition}
Up to a scalar function the expression \eqref{eq:psi} for $\psi$ can be directly
obtained, as follows:
\begin{exercise}
Over the set where $\Ad_g+1$ is invertible, the vector fields 
$\f{\xi^L+\xi^R}{2}$ span the tangent space. Hence, there is a unique
2-form $\varsigma$ on this set, with
\[ \iota(\ts{\f{\xi^L+\xi^R}{2}})\varsigma+B(\ts{\f{\theta^L-\theta^R}{4}},\xi)
=0.\]
Deduce that $e^\varsigma$ is a pure spinor defining $F$, hence
coincides with $\psi$ up to a scalar function. Next, check that 
\[\varsigma=\ts{\f{1}{4} }B\Big(\ts{\f{1-\Ad_g}{1+\Ad_g}}\theta^L,\theta^L\Big)\]
is the unique solution of the defining equation for $\varsigma$. 
\end{exercise}

\vskip.5in
\section{Dirac structures}\label{sec:dirac}
\subsection{Courant's integrability condition}
One of Courant's main discoveries in \cite{cou:di} was the existence
of a natural integrability condition for almost Dirac structures
$E\subset \T M$.  
Following Alekseev-Xu \cite{al:der} and Gualtieri \cite{gua:ge}, we will express the
Courant integrability condition in terms of the spinor representation.
The Lagrangian sub-bundle $E$ defines a filtration on the spinor module $\Omega(M)$:
\[ \Omega(M)=\Omega^{(n)}(M)\supset\cdots \Omega^{(1)}(M)\supset
\Omega^{(0)}(M).\]
Here $\Omega^{(k)}(M)$ consists of differential forms $\gamma$ with
$\rho(w_0)\cdots\rho(w_k)\gamma=0$ for all $w_i\in\Gamma(E)$.

Let us fix a closed 3-form $\eta\in\Omega^3(M)$ (possibly zero). Note
that $\d+\eta$ is again a differential.
\begin{lemma}
Let $\phi\in\Omega(M)$ be a (locally defined) pure spinor with $N_\phi=E$. Then
\[(\d+\eta)\phi\in \Omega^{(3)}(M).\]  
\end{lemma}
\begin{proof}
Let $w_i\in\Gamma(E)$. Since $\rho(w_i)$ annihilates $\phi$,
we have 
\begin{equation}\label{eq:3} 
\rho(w_1)\rho(w_2)\rho(w_3)(\d+\eta) \phi=[\rho(w_1),[\rho(w_2),[\rho(w_3),\d+\eta]]]\phi,
\end{equation}
using graded commutators of operators on $\Omega(M)$. 
A calculation (cf. Exercise \ref{ex:1} below) shows that the triple
commutator of operators is multiplication by a smooth function. 
Thus $\rho(w_1)\rho(w_2)\rho(w_3)(\d+\eta) \phi$ is a function times $\phi$,
and hence is annihilated by $\rho(w_0)$.
\end{proof}

%\[ Y(w_1,w_2,w_3)=\iota(v_1)\iota(v_2)(\d\alpha_3)-\iota(v_1)L(v_3)\alpha_2
%+\iota([v_2,v_3])\alpha_1+\iota(v_1)\iota(v_2)
%\iota(v_3)\eta.\]
%
%Here we have written $w_i=v_i\oplus\alpha_i$. 

We may now state the Courant integrability condition. 
\begin{definition}
An almost Dirac structure $E\subset \T M$ is called integrable
relative to the closed 3-form $\eta$ if, for
any (locally defined) pure spinor $\phi$ with $N_\phi=E$, 
\begin{equation}\label{eq:courint}
 \on{gr}^3\big((\d+\eta)\phi)=0.\end{equation}
\end{definition}
Note that this condition does not depend on the choice of $\phi$,
since 
\[ \on{gr}^3\big((\d+\eta)(f\phi))=f\on{gr}^3((\d+\eta)(\phi).\]
By \eqref{eq:courint}, $E$ is integrable
if and only if $(\d+\eta)\phi\in\Omega^{(2)}(M)$. Since $\phi$ and
$(\d+\eta)\phi$ have opposite parity, this is in fact equivalent to
the condition
\begin{equation} (\d+\eta)\phi\in\Omega^{(1)}(M).\end{equation}
\begin{definition}
A \emph{Dirac manifold} is a triple $(M,E_M,\eta_M)$, consisting of a
manifold $M$, an almost Dirac structure $E_M$, and a closed 3-form
$\eta_M$ such that $E_M$ is integrable relative to $\eta_M$. A smooth
map $\Phi\colon M\to M'$ between two Dirac manifolds is called a
\emph{(strong) Dirac map} if each $\d_x\Phi\colon T_xM\to T_{\Phi(x)}M'$ is a
linear (strong) Dirac map, and in addition 
\[ \Phi^*\eta_{M'}=\eta_M.\]
\end{definition}

\begin{remark}
The integrability condition may be rephrased as
      $(\d+\eta)\phi=\rho(w)\phi$ for \emph{some} section
      $w\in\Gamma(\T M)$. It is not always possible to choose $\phi$
      in such a way that $(\d+\eta)\phi=0$. As shown by Alekseev-Xu
      \cite{al:der}, the obstruction is the `modular class' of $E$.
\end{remark}

\subsection{Examples}
\begin{examples}\label{ex:exam}
\begin{enumerate}
\item
Let $\omega$ be a 2-form and $\phi=e^{-\omega}$. Then $(\d+\eta)\phi=(-\d\omega+\eta)\wedge
\phi$ lies in $\Omega^{(1)}(M)$ if and only if $\d\omega=\eta$. 
From now on, we will view any manifold $M$ with 2-form $\omega$ as a Dirac
manifold, taking $E_M=\on{Gr}_\omega$. 
Observe that $\Phi\colon M\to \pt$ is a strong
Dirac map if and only if $\omega$ is symplectic (closed and
non-degenerate). 
\item More generally, if $E$ is integrable with respect to $\eta$, and
      $\tau$ is any 2-form, then $E^\tau$ is integrable with respect
      to $\eta+\d\tau$.
\item Let $\pi$ be a bi-vector field and $\mu_M$ a volume form on $M$.
      Then $\phi=e^{-\iota(\pi)}\mu_M$ satisfies 
\[
(\d+\eta)\phi=\iota\big(-\hh[\pi,\pi]_{\on{Sch}}-\pi^\sharp(\eta)+X_\pi+Y_{\pi,\eta})\phi.\]
Here $[\cdot,\cdot]_{\on{Sch}}$ is the Schouten bracket on
multi-vector fields, $X_\pi\in\X^1(M)$ is the vector field 
defined by $\d\iota(\pi)\mu_M=-\iota(X_\pi)\mu_M$, 
$\pi^\sharp$ is
the bundle map from $\wedge T^*M$ to $\wedge TM$, 
and $Y_{\pi,\eta}$ is
the vector field $Y_{\pi,\eta}=\pi^\sharp(\iota(\pi)\eta)$. 
The Courant integrability condition reduces to the condition
\[\hh[\pi,\pi]_{\on{Sch}}+\pi^\sharp(\eta)=0,\]
defining a \emph{twisted Poisson structure}. These structures were
introduced by Klimcik-Strobl \cite{kli:wzw} and further studied by
\v{S}evera-Weinstein \cite{sev:po}. It was argued by
Kosmann-Schwarzbach-Laurent-Gengoux \cite[Theorem 6.1]{kos:mod} (see
also \cite[Example 6.2]{al:der}) that the sum $X_\pi+Y_{\pi,\eta}$
plays the role of the \emph{modular vector field} for a twisted
Poisson structure.

%
%From now on, any Poisson
%manifold $(M,\pi)$ (i.e. $[\pi,\pi]_{\on{Sch}}=0$) will be viewed as a
%Dirac manifold with $E=\on{Gr}_\pi$ and $\eta_M=0$.  Note that any
%Poisson map is a strong Dirac map.
%
\item Take $\eta=0$, and let
      $\alpha_1,\ldots,\alpha_{k}\in\Omega^1(M)$ be a collection of
      pointwise linearly independent 1-forms, and $K\subset TM$ be the
      codimension $k$ distribution given as the intersection of their
      kernels. Then $\phi=\alpha_1\wedge\cdots \wedge \alpha_{k}$ is
      pure spinor defining $E=K\oplus \on{ann}(K)$. The integrability
      condition $\d\phi\in\Omega^{(1)}(M)$ holds if and only if
      $\d\phi=\beta\wedge\phi$ for a 1-form $\beta$.  This is one
      version of the standard (Frobenius) integrability condition for
      distributions.
\item \label{it:e} The Courant integrability condition has an obvious
      generalization to complex almost Dirac structures $E\subset \T
      M\otimes\C$. Given an almost complex structure on $M$, i.e.  a
      linear complex structure on the tangent bundle, $J\in
      \Gamma(\on{End}(TM))$, $J^2=-\on{Id}_{TM}$, one obtains a 
      linear complex structure $\mathbb{J}=J\oplus (-J^*)\in \Gamma(\on{End}(\TM)),\ 
      \ \mathbb{J}^2=-\on{Id}_{\TM}$.  Let $E\subset \TM\otimes\C$ be
      the $+i$ eigenbundle of $\mathbb{J}$.  It turns out that $E$ is
      Courant integrable if and only if the almost complex structure
      $J$ is integrable, i.e comes from complex coordinate charts with
      holomorphic transition functions.  This is the motivating
      example for the \emph{generalized complex geometry}, developed
      by Hitchin \cite{hi:gen} and Gualtieri \cite{gua:ge}.
\end{enumerate}
\end{examples}

\begin{exercise}
In Example \ref{ex:exam}\eqref{it:e}, give a pure spinor
$\phi\in\Om(M)\otimes\C$ defining $E$.
\end{exercise}
\begin{exercise}
(See \cite{ev:po}.) Work out a formula for 
\[ e^{\iota(\pi)}\circ \d \circ e^{-\iota(\pi)}=\d+[\iota(\pi),\d]+
\hh [\iota(\pi),[\iota(\pi),\d]]+\cdots.\]
(The resulting expression contains terms at most quadratic in $\pi$.) 
Use this to show 
\[ \d(e^{-\iota(\pi)}\mu_M)=\iota(-\hh [\pi,\pi]_{\on{Sch}}+X_\pi)\mu_M\]
for any volume form $\mu_M$. Similarly show that  
\[ \eta \wedge (e^{-\iota(\pi)}\mu_M)=
\iota\big(\!-\pi^\sharp(\eta)+Y_{\pi,\eta}\big) (e^{-\iota(\pi)}\mu_M).\]
\end{exercise}

\begin{exercise}\label{ex:1}
\begin{enumerate}
\item 
Verify that the following formula defines a bilinear map $\Cour{\cdot,\cdot}\colon
\Gamma(\T M)\times\Gamma(\T M)\to \Gamma(\TM)$: 
\[ \rho(\Cour{w_1,w_2})=[\rho(w_1),[\rho(w_2),\d+\eta]].\]
This is the definition of the (non skew-symmetric)
\emph{Courant bracket} $\Cour{\cdot,\cdot}$ on $\Gamma(\TM)$ as a \emph{derived bracket}.
See Roytenberg \cite{roy:co}, Alekseev-Xu \cite{al:der} and
Kosmann-Schwarzbach \cite{kos:der}.

\item
Conclude that for any $w_1,w_2,w_3\in\Gamma(\T
M)$, the operator 
\[ [\rho(w_1),[\rho(w_2),[\rho(w_3),\d+\eta]]]\]
on $\Omega(M)$ is multiplication by the smooth function
$Y(w_1,w_2,w_3)$. 
\item 
Show that for any almost Dirac structure $E\subset \T M$, the
restriction of $Y$ to sections of $E$ defines an anti-symmetric
tensor $Y_E\in\wedge^3 E^*$.
\end{enumerate} 
\end{exercise}

\begin{proposition}\label{prop:liealgebroid}
The almost Dirac structure $E$ is integrable if and only if
$\Gamma(E)$ is closed under Courant bracket $\Cour{\cdot,\cdot}$. 
In this case, the restriction $[\cdot,\cdot]_E$ of the Courant bracket
to $\Gamma(E)$ 
defines a Lie algebroid structure on $E$: That is, it is a Lie
bracket, the projection map $a\colon\Gamma(E)\to \mf{X}(M)$ is a Lie algebra
homomorphism, and 
\[ [w_1,fw_2]_E=f[w_1,w_2]_E+v_1(f)\,w_2,\ \  w_i\in\Gamma(E)
\]
where $v_1=a(w_1)$. 
\end{proposition}
\begin{proof}
Since $E$ is Lagrangian, we have $\Cour{w_2,w_3}\in\Gamma(E)$ for
all $w_2,w_3\in\Gamma(E)$ if and only if 
\[Y(w_1,w_2,w_3)= [\rho(w_1),\rho(\Cour{w_2,w_3})]
=\l w_1,\Cour{w_2,w_3}\r=0\]
for all $w_1,w_2,w_3\in\Gamma(E)$. The remaining claims are left as an exercise.
\end{proof}
The theory of Lie algebroids \cite{ca:ge,mac:gen} shows that the
generalized distribution $\on{ran}(E)=\pr_{TM}(E)$ is
\emph{integrable}, i.e. defines a generalized foliation. Moreover, the
leaves $S\subset M$ of this foliation carry \emph{2-forms}
$\omega_S\in\Omega^2(S)$, defined pointwise by \eqref{eq:twoform}.

For $E=\on{Gr}_\pi$ the graph of a Poisson bi-vector field (i.e.
$\eta=0$), this is just the usual foliation by symplectic leaves
$S\subset M$, with $\omega_S$ the symplectic 2-forms. More generally,
in the twisted Poisson case $\hh [\pi,\pi]_{\on{Sch}}+\pi^\sharp\eta=0$ one
still obtains a foliation. The 2-forms on the leaves are again
non-degenerate (since $E\cap TM=\{0\}$), but are not closed 
in general:
\begin{proposition}\label{prop:tform}
  Let $E\subset \TM$ be a Dirac structure (relative to a closed 3-form
  $\eta\in\Omega^3(M)$). Then the 2-forms
  $\omega_S$ on the leaves $S\subset M$ satisfy
  $\d\omega_S=\iota_S^*\eta$.
\end{proposition}
\begin{proof}
  Given any point $x\in S$, we may pass to a neighbourhood of $x$ to
  reduce to the case $M=S\times N$, with $\iota_S$ the inclusion as
  $S\times\{y\}$ for some $y\in N$. View $\omega_S$ as a form on
  $S\times N$, and define $\gamma:=e^{\omega_S}\phi$.
  Then $\gamma|_S$ is a nowhere vanishing section of the top exterior
  power of $T^*N|_S\cong \on{ann}(S)\subset T^*M|_S$.  By assumption,
  there exists a vector field $v$ and a 1-form $\alpha$ such that
  $(\d+\eta)\phi=\iota(v)\phi+\alpha\phi$. This yields:
\[0=e^{\omega_S}(\d+\eta-\iota(v)-\alpha)\phi
   =\big(\eta-\d\omega_S+\iota(v)\omega-\alpha\big)\gamma
   +(\d\gamma-\iota(v)\gamma).\]
Restricting to $S$, and taking the component in
$\Gamma(\wedge^3 T^*S\otimes \wedge^{\on{top}}T^*N|_S)$ we find 
$(\iota_S^*\eta-\d\omega_S)\gamma|_S=0$. Hence
$\iota_S^*\eta=\d\omega_S$. 
\end{proof}

\subsection{Integrability of the Cartan-Dirac structure}
Let us now return to the example of a Lie group $G$ with an invariant inner
product $B$ on $\g$. Suppose that $G$ admits an invariant
orientation and that $\Ad\colon G\to \on{O}(\g)$ lifts to the
$\on{Pin}$ group. Let $\phi,\psi\in\Omega(G)$ be the pure spinors 
defining the almost Dirac structures $E,F$. By construction, both 
$\phi$ and $\psi$ are $\Ad$-invariant differential forms. 

Now let $\eta\in\Omega^3(G)$ be the left-invariant 3-form
\[ \eta=\f{1}{12} B(\theta^L,[\theta^L,\theta^L]).\]
Since $B$ is invariant, one may replace $\theta^L$ with $\theta^R$ in
this formula, thus $\eta$ is also right-invariant. In particular,
$\eta$ is closed (since any bi-invariant differential form on a Lie
group is closed). Letting $\xi^\sharp=\xi^L-\xi^R$ be the generating
vector fields for the conjugation action, one finds 
\[ \iota(\xi^\sharp)\eta=-\d\ B(\ts{\f{\theta^L+\theta^R}{2}},\xi)
\]
As a consequence, we see that the commutator of $\d+\eta$ with the 
generating sections $e(\xi)$ of $E$ are, 
\[[\rho(e(\xi)),\d+\eta]=[\iota(\xi^\sharp)+B(\ts{\f{\theta^L+\theta^R}{2}},\xi),\d+\eta]
=L(\xi^\sharp).\]
(Here $L(X)=[\iota(X),\d]$ denotes the Lie derivative in the direction of a vector
field $X$.) It hence follows that
\[ \rho(e(\xi))(\d+\eta)\phi=[\rho(e(\xi)),\d+\eta]\phi=L(\xi^\sharp)\phi=0.\]
Thus $(\d+\eta)\phi\in\Omega^{(0)}(M)$. Since the parity of
$(\d+\eta)\phi$ is opposite to that of $\phi$, we obtain:
\begin{theorem}
The pure spinor $\phi$ satisfies 
\[ (\d+\eta)\phi=0.\]
\end{theorem}
In particular, we see that $E$ is a Dirac structure.
\begin{definition}
The Dirac structure $E$ on $G$ is called the \emph{Cartan-Dirac structure}.
\end{definition}
The integrability of $E$ explains our earlier observation that the
distribution $\on{ran}(E)$ is just the tangent distribution for the
generalized foliation by conjugacy classes. Furthermore, Proposition
\ref{prop:tform} tells us that the GHJW 2-form $\omega_\Co$ on the
conjugacy classes satisfies,
\[ \d\omega_\Co=\iota_\Co^*\eta.\]
\begin{remark}
The Cartan-Dirac structure was discovered independently by Anton
Alekseev, Pavol \v{S}evera and Thomas Strobl, around the end of the last
century.
\end{remark}

\begin{remark}
By contrast, the almost Dirac structure $F$ is not
integrable. Instead, one has \cite{al:pur}
\[ (\d+\eta)\psi=\rho(e(\Xi))\psi\]
where $\Xi\in\wedge^3\g$ is the `structure constants tensor', and
$e(\Xi)\in\Gamma(\wedge^3 E)$ is defined using the extension of
$e\colon \g\to \Gamma(E)$ to an algebra homomorphism
$\wedge\g\to\Gamma(\wedge(E))$.
\end{remark}

\vskip.5in
\section{Group-valued moment maps}
The theory of $G$-valued moment maps was introduced in the paper
\cite{al:mom}.  One of its main applications is that it provides a
natural framework for the construction of symplectic forms on moduli
spaces of flat connections. 

\subsection{Definition of q-Hamiltonian $G$-spaces}
Let $G$ be a connected Lie group with a bi-invariant pseudo-Riemannian 
metric, and let $B$ be the corresponding invariant inner product on $\g$. 
Let $M$ be a manifold. A \emph{$G$-action} on $M$ is a group
homomorphism $\A\colon G\to\on{Diff}(M)$ such that the action map
$G\times M\to M,\ (g,x)\mapsto \A(g).x$ is smooth. Similarly, a
\emph{$\g$-action} is a Lie algebra homomorphism $\A\colon \g\to
\mf{X}(M)$ such that the map $\g\times M\to TM,\ (\xi,x)\mapsto
\A(\xi)_x$ is smooth. We will write $\xi^\sharp=\A(\xi)$. 
For any $G$-action, the generating vector fields
(defined with the appropriate sign) give a $\g$-action. Conversely, if
$M$ is compact and $G$ is simply connected, any $\g$-action integrates
to a $G$-action.

\begin{definition}\cite{al:mom} \label{def:def1}
A Hamiltonian $\g$-space with $G$-valued moment map is a $\g$-manifold $M$,
together with a $\g$-invariant 2-form $\omega\in\Omega^2(M)$ and a 
$\g$-equivariant map $\Phi\in C^\infty(M,G)$ such that
\begin{enumerate}
\item $\d\omega=\Phi^*\eta$, 
\item $\iota(\xi^\sharp)\omega=\Phi^* B(\ts{\f{\theta^L+\theta^R}{2}},\xi),\ \ \ \ \xi\in\g$
\ \ (Moment map condition.)  
\item $\ker(\omega_x)=\{\xi^\sharp(x)|\ \Ad_{\Phi(x)}\xi=-\xi\},\ x\in M$\ \ (Minimal degeneracy condition.)
\end{enumerate}
\end{definition}

\begin{remark}
  As pointed out in \cite{al:mom}, (b) is the simplest $G$-valued
  analogue to the defining property for $\g^*$-valued moment maps
  $\Phi_0\colon M\to \g^*$,
  $\iota(\xi^\sharp)\omega_0=-\d\l\Phi_0,\xi\r$.  It follows from the
  work of Bursztyn-Crainic \cite{bur:di} and Xu \cite{xu:mom} (see
  also \cite{al:pur}) that (c) may be replaced by
  the more elegant condition,
\[ \ker(\omega_x)\cap \ker(\d_x\Phi)=0.\]
\end{remark}

The theory of $G$-valued moment maps was developed in \cite{al:mom},
and subsequent papers, in full analogy to the familiar theory of
$\g^*$-valued moment maps. However, the proofs were much more
complicated than in the $\g^*$-valued theory, and for technical
reasons it was necessary to assume that $B$ is positive definite.
Unfortunately, this restriction excludes several interesting examples,
such as representation varieties for non-compact semi-simple Lie
groups. (The Killing form of such groups is indefinite.)  In the
following approach to group-valued moment maps via Dirac structures
these difficulties are no longer present.
\begin{theorem}[Bursztyn-Crainic] Definition \ref{def:def1} is equivalent to the following Definition \ref{def:def2}.
\end{theorem}

\begin{definition}\label{def:def2} A Hamiltonian $\g$-space with $G$-valued moment map is a manifold 
$M$ with a 2-form $\omega$, together with a strong Dirac map
$\Phi\colon M\to G$. 
\end{definition}
Here $M$ is viewed as a Dirac manifold with $E_M=\on{Gr}_\omega$ and
3-form $\eta_M=\d\omega$, while $G$ carries the Cartan-Dirac
structure.  Recall that $\Phi^*\eta=\eta_M$ as part of the definition
of a Dirac map from $(M,\on{Gr}_\omega,\d\omega)$ to $(G,E,\eta)$.

Note that Definition \ref{def:def2} no longer
mentions the $\g$-action on $M$, the equivariance of $\omega$ and
$\Phi$, or the minimal degeneracy property: as shown by
Bursztyn-Crainic, all of this comes for free!

\begin{remarks}
One is immediately led to consider arbitrary Dirac
manifolds $(M,E_M,\eta_M)$ together with strong Dirac maps $M\to G$. 
As shown by Bursztyn-Crainic, one recovers the
theory of q-Poisson manifolds \cite{al:ma,al:qu}. Definition
\ref{def:def2} is parallel to the definition of Hamiltonian
$\g$-spaces with $\g^*$-valued moment maps: These may be defined as
manifolds $M$ with closed 2-forms $\omega$ and strong Dirac maps
$\Phi\colon M\to \g^*$. Here $\g^*$ carries the Dirac structure coming
from its Kirillov-Poisson structure. Similarly, Lu's notion
\cite{lu:mo} of moment maps $\Phi\colon M\to G^*$ for Poisson
$G$-actions on symplectic manifolds (where $G,G^*$ are dual Poisson
Lie groups) can be phrased in this way.
\end{remarks}

In most cases of interest, the $\g$-action on $M$ exponentiates to an action of
$G$: 
\begin{definition}
  Let $M$ be a $G$-manifold, together with a $G$-invariant 2-form
  $\omega$ and a $G$-equivariant map $\Phi\colon M\to G$. Then
  $(M,\omega,\Phi)$ is called a Hamiltonian $G$-space with $G$-valued
  moment map, or simply a \emph{q-Hamiltonian $G$-space}, if $\Phi$ is
  a Dirac map, and the $\g$-action generated by $\Phi$ is the
  infinitesimal $G$-action. 
\end{definition}

\begin{example}
Every conjugacy class $\Co\subset G$, equipped with the GHJW 2-form, 
is a q-Hamiltonian $G$-space, with moment map the inclusion. 
\end{example}

\subsection{Volume forms}
Definition \ref{def:def2} greatly simplifies many of the constructions 
with $G$-valued moment maps. For instance, generalizing our arguments 
for the {\Fact}  about conjugacy classes, one obtains 
the following 
\begin{theorem}\cite{al:pur}
Assume that the homomorphism $\Ad\colon G\to \on{O}(\g)$ lifts to the
group $\on{Pin}(\g)$, and let $\psi\in\Omega(G)$ be the pure spinor 
with $N_\psi=F$, defined by this lift $\wt{\Ad}$. 
Let $(M,\omega,\Phi)$ be a q-Hamiltonian $G$-space. Then
\[ (e^\omega\Phi^*\psi)_{[\on{top}]}\in \Omega(M)\]
is a $G$-invariant volume form. 
\end{theorem}
If $G$ is connected, we see that $\dim M$ must be even (since $\psi$
is an even form in this case). 

\subsection{Products}
Let us next consider products of q-Hamiltonian $G$-spaces. 

For ordinary Hamiltonian $G$-spaces $(M_i,\omega_i)$ with moment
maps $\Phi_i\colon M_i\to \g^*$, the product is simply the direct
product $M_1\times M_2$ with the diagonal action, the sum of the 2-forms
$\omega_1+\omega_2$ and the sum $\Phi_1+\Phi_2$ of the moment maps.
Similarly, if $G$ is a Poisson Lie group with dual Poisson-Lie group
$G^*$, the product operation for Lu's Hamiltonian $G$-spaces with
$G^*$-valued moment maps takes the sum of the 2-forms and the
pointwise product of the moment maps. In this case, the 
$G$-action is a certain twist of the diagonal action, see 
\cite{fl:pc,lu:mo}.

These two constructions work because the addition map
$\on{Add}\colon\g^*\times\g^* \to \g^*$, respectively the product map
$\Mult\colon G^*\times G^*\to G^*$, are Poisson. For $G$-valued moment
maps, the situation is slightly different since the multiplication map
$\Mult\colon G\times G\to G$, \emph{as it stands}, is not a strong Dirac map
if one simply takes the direct product Dirac structure on $G\times G$.
Instead, the product operation involves a \emph{gauge transformation}. 

\begin{definition}
Let $E_M$ be an almost Dirac structure on $M$, defined (locally) by a
pure spinor $\phi_M$, and $\tau\in \Omega^2(M)$ a 2-form. Then the
\emph{gauge transformation} $E_M^\tau$ is the almost Dirac structure
defined (locally) by the pure spinor $e^{-\tau}\phi_M$.
\end{definition}
Note that if $E_M$ is integrable with respect to a closed 3-form 
$\eta_M$, then $E_M^\tau$ is integrable with respect to $\eta_M+\d\tau$. 
In our case, we need a suitable gauge transformation of 
$E_{G\times G}:=E_G\times E_G$. Let 
\[ \tau:=\hh B(\pr_1^*\theta^L,\pr_2^*\theta^R)\in\Omega^2(G\times G).\]
This 2-form has the property \cite{we:sy}, 
\[ \Mult^*\eta=\pr_1^*\eta+\pr_2^*\eta+\d\tau.\]
\begin{theorem}\cite{al:pur}\label{th:product}
The multiplication map $\on{Mult}\colon G\times G\to G$ is a strong
Dirac map from $(G\times G,E_{G\times G}^\tau)$ to $(G,E_G)$.
\end{theorem}
One may use this result to define the \emph{fusion product} of two
q-Hamiltonian $G$-spaces $M_1,M_2$, or more generally to pass to the
diagonal action in a q-Hamiltonian $G\times G$-space $M$ (e.g.
$M=M_1\times M_2$). Indeed let $(M,\omega,\Phi)$ be a q-Hamiltonian
$G\times G$-space with moment map $\Phi=\Phi_1\times\Phi_2$. Put
\[ \Phi^{\on{fus}}=\Phi_1\Phi_2,\  
\omega^{\on{fus}}=\omega+(\Phi_1,\Phi_2)^*\tau.\]
Then $(M,\omega^{\on{fus}},\Phi^{\on{fus}})$, with diagonal
$G$-action, is a q-Hamiltonian $G$-space. This follows rather easily
from Theorem \ref{th:product}, since the composition of two strong
Dirac maps is again a strong Dirac map.

\begin{remark}
  For the case of compact Lie groups, and working with the Definition
  \ref{def:def1}, this result was obtained in \cite{al:mom} by a
  fairly complicated argument. The main difficulty in this approach
  was to show that $\omega^{\on{fus}}$ is again minimally degenerate:
  It is not easy to compute the kernel of $\omega^{\on{fus}}$ by
  `direct calculation'!
\end{remark}

\subsection{Exponentials}
Let the dual of the Lie algebra $\g^*$ be equipped with the
Kirillov-Poisson structure $\pi$. Its graph $\on{Gr}_\pi$ defines a
Dirac structure. Use the inner product $B$ to identify 
$\g^*\cong\g$. Let $\varpi\in\Omega^2(\g)$ be the 2-form, 
obtained by applying the de Rham homotopy operator to 
$\exp^*\eta$. Thus $(\on{Gr}_\pi)^\varpi$ is a Dirac structure 
relative to the closed 3-form $\d\varpi=\exp^*\eta$. Now let 
$\g_\natural \subset \g$ be the open subset where 
$\exp$ is a local diffeomorphism. 
\begin{theorem}\cite{al:pur}
The restriction of $\exp$ to the subset $\g_\natural$ is a strong Dirac map,
relative to the Dirac structures $(\on{Gr}_\pi)^\varpi$ on $\g$ and
the Cartan-Dirac structure on $G$.
\end{theorem}

Suppose now that $(M,\omega_0,\Phi_0)$ is an ordinary Hamiltonian
$G$-space (thus $\omega_0$ is a symplectic form, and $\Phi_0\colon
M\to \g^*\cong\g$ a moment map in the usual sense). Let
$\omega=\omega_0+\Phi_0^*\varpi$ and $\Phi=\exp\circ \Phi_0$.  Then
$(M,\omega,\Phi)$ is a q-Hamiltonian $G$-space \emph{provided that}
$\Phi_0(M)\subset \g_\natural$. Again, this just follows from the fact
that the composition of two strong Dirac maps is again a strong Dirac
map. Conversely, suppose $U\subset \g_\natural$ is an open subset
where $\exp$ is a diffeomorphism (with inverse denoted $\log$), and
$(M,\omega,\Phi)$ is a q-Hamiltonian $G$-space. Put
$\Phi_0=\log(\Phi)$ and $\omega_0=\omega-\Phi_0^*\varpi$. Then
$\omega_0$ is symplectic, and $(M,\omega_0,\Phi_0)$ is a Hamiltonian
$G$-space in the usual sense. For $(M,\omega_0,\Phi_0)$ one has all
the standard results from symplectic geometry, which one may then
translate back to the q-Hamiltonian setting. For instance, 
the Meyer-Marsden-Weinstein reduction theorem for Hamiltonian manifolds
(see Bates-Lerman \cite{bat:pro} for a very general version) 
yields:
\begin{proposition}[Symplectic reduction of q-Hamiltonian manifolds]
Let $(M,\omega,\Phi)$ be a q-Hamiltonian $G$-space, with proper moment map. 
Suppose the action of $G$ is proper, and that $e$ is a regular value of the moment map. Then $G$ acts locally free on $\Phinv(e)$, and the reduced space 
\[ M\qu G=\Phinv(e)/G\]
is a symplectic orbifold. (If $e$ is not a regular value, $M\qu G$ 
is a stratified symplectic space.) 
\end{proposition}

\subsection{Examples}
\subsubsection{Homogeneous spaces}\label{subsec:sym}
Let $G$ be a Lie group with involution $\sig\in\on{Aut}(G)$,
$\sig^2=1$, and consider the symmetric space $M=G/G^\sigma$.  Let
$\wh{G}=\Z_2\ltimes G$ be the semi-direct product, defined using the
action of $\Z_2=\{1,\sigma\}$ on $G$. Then $M$ may be viewed as the
conjugacy class of the element $(\sigma,e)\in\wh{G}$. Hence, if $\g$
carries an invariant scalar product which is preserved under the
involution, the space $M$ becomes a q-Hamiltonian 
$\wh{G}$-space, with moment map
\[ \Phi\colon M\to \wh{G},\ \ gG^\sigma\mapsto (\sig,g\sig(g)^{-1}).\]
Since $(\sig,e)$ squares to the group unit, Exercise \ref{ex:square}
shows that the 2-form $\omega$ on $M$ is identically zero. Note also
that the action of $\Z_2\subset \wh{G}$ is trivial, so that the action
of $\wh{G}$ descends to $G$.

Consider now the fusion product of $M$ with itself. Letting
$\Phi_i=\Phi\circ \pr_i\colon M\times M\to \wh{G}$, the map
$\Phi=\Phi_1\Phi_2$ takes values in the subgroup $G\subset \wh{G}$:
\[ \Phi(g_1G^\sigma,\ g_2G^\sigma)=\sig(g_1)g_1^{-1} g_2 \sig(g_2^{-1}).\]
Hence $M\times M$ is a q-Hamiltonian $G$-space. 
\subsubsection{The double}\label{ex:g}
Any Lie group $G$ may be viewed as a symmetric space for the group
$G\times G$, with action $(g_1,g_2).a=g_1 a g_2^{-1}$. Here
$\sigma\in\on{Aut}(G\times G)$ is the involutive automorphism
$\sig(g_1,g_2)=(g_2,g_1)$, fixing the diagonal subgroup, and the
inclusion of the first factor identifies the quotient $(G\times G)/G$
with $G$. Hence, given an invariant scalar product on $\g$, the
example in Section \ref{subsec:sym} shows that $G$ is a q-Hamiltonian $\Z_2\ltimes(G\times
G)$-space. Taking a fusion product of $G$ with itself, we find that
$D(G):=G\times G$ is a q-Hamiltonian $G\times G$-space, with action
\[ (g_1,g_2).(a,b)=(g_1 a g_2^{-1},\ g_2 b g_1^{-1})\]
and moment map $(a,b)\mapsto (ab,a^{-1}b^{-1})$. The space $D(G)$ is
called the \emph{double} of $G$. 
\begin{remark}
  The double $D(G)$ is the counterpart, in the q-Hamiltonian category,
  of the cotangent bundle $T^*G$ in the usual Hamiltonian category. In
  fact, as observed in Bursztyn-Crainic-Weinstein-Zhu \cite{bur:int}
  the double $D(G)\rightrightarrows G$ (viewed as a groupoid over $G$,
  with source and target maps the two components of the moment map)
  `integrates' the Dirac manifold $G$ in a similar sense as
  $T^*G\rightrightarrows \g^*$ integrates the Poisson manifold $\g^*$.
  Ping Xu \cite{xu:mom} presents $D(G)\rightrightarrows G $ as an
  example of a \emph{quasi-symplectic groupoid}.
\end{remark}
Passing to the diagonal action, $G\times G$
becomes a q-Hamiltonian $G$-space with moment map the group
commutator: 
\[ (a,b)\mapsto aba^{-1}b^{-1}.\]
This is called the \emph{fused double}, denoted  $\ti{D}(G)$. 
Taking a fusion product of several copies of $\ti{D}(G)$ with itself, 
the space $G^{2h}$ becomes a q-Hamiltonian $G$-space with moment map 
\[ \Phi\colon (a_1,b_1,\ldots,a_h,b_h)\mapsto \prod_{i=1}^h a_i b_i a_i^{-1} b_i^{-1}.\]
The symplectic quotient $M\qu G$ is just the representation variety for a 
closed oriented surface of genus $h$:
\[ M\qu G=\on{Hom}(\pi_1(\Sig),G)/G\]
Equivalently, $M\qu G$ is the moduli space of flat principal
$G$-bundles on $\Sig$. It was shown in \cite{al:mom} that the
symplectic structure obtained by this finite-dimensional reduction,
coincides with that coming from Atiyah-Bott's \cite{at:ge} gauge
theory construction. 
More generally, if $\Co_1,\ldots,\Co_r$ are conjugacy
classes in $G$, the symplectic quotient 
\[ (G^{2s}\times\Co_1\times\cdots\times\Co_r)\qu G\]
is the moduli space of flat $G$-bundles over an oriented surface
$\Sig$ of genus $h$ with $r$ boundary components, with restrictions to
the $j$th boundary component $(\partial\Sig)_j\cong S^1$ in the given
conjugacy classes. (Note $\Hom(\pi_1(S^1),G)/G)=G/\Ad(G)$ is the set
of conjugacy classes.) 

\begin{remark}
We stress that no compactness assumption is needed for these results. 
In fact, one could even work over the complex numbers, and obtain a 
complex symplectic structure over the representation variety for a 
complex Lie group. 
\end{remark}

\subsubsection{Spheres}
There are other examples of q-Hamiltonian spaces which are unrelated
to moduli spaces, such as various examples of \emph{multiplicity-free}
q-Hamiltonian spaces. Let $\SU(n)$ act on $\C^n$ in the standard way,
and consider the unit sphere $S^{2n}\subset \C^n\times\R$ with the
restricted action.  By a result of Hurtubise-Jeffrey-Sjamaar
\cite{hu:imp}, there exists an invariant 2-form $\omega$ and an
equivariant map $\Phi\colon S^{2n}\to \SU(n)$ for which
$(S^{2n},\omega,\Phi)$ is a q-Hamiltonian $\SU(n)$-space. (The special
case $n=2$ was discussed in \cite{al:du}.) 
\vskip1in

\bibliographystyle{amsplain}
%\bibliography{../Biblio/ref}
%\bibliography{../../Biblio/ref}

\def\cprime{$'$} \def\polhk#1{\setbox0=\hbox{#1}{\ooalign{\hidewidth
  \lower1.5ex\hbox{`}\hidewidth\crcr\unhbox0}}} \def\cprime{$'$}
  \def\cprime{$'$} \def\polhk#1{\setbox0=\hbox{#1}{\ooalign{\hidewidth
  \lower1.5ex\hbox{`}\hidewidth\crcr\unhbox0}}} \def\cprime{$'$}
  \def\cprime{$'$}
\providecommand{\bysame}{\leavevmode\hbox to3em{\hrulefill}\thinspace}
\providecommand{\MR}{\relax\ifhmode\unskip\space\fi MR }
\providecommand{\MRhref}[2]{%
  \href{http://www.ams.org/mathscinet-getitem?mr=#1}{#2}
}
\providecommand{\href}[2]{#2}

\end{document}